\documentclass[a4paper,12pt]{article}
\usepackage{amsmath,amsthm,verbatim,amssymb,amsfonts,amscd, graphicx}
\usepackage[all]{xy}
\usepackage{xcolor}
\usepackage{geometry}
\usepackage{graphicx}
\usepackage{hyperref}
\usepackage{enumitem}
\usepackage{indentfirst}
\definecolor{darkblue}{rgb}{0,0,0.3}
\definecolor{urlblue}{rgb}{0,0,0.7}
\hypersetup{
	colorlinks=true,
	citecolor=blue,
	linkcolor=darkblue,
	urlcolor=urlblue,
}

\newcommand{\RR}{\mathbb{R}}

\newcommand{\NN}{\mathbb{N}}

\let\div\undefined

\DeclareMathOperator{\diam}{diam}
\DeclareMathOperator{\div}{div}

\DeclareMathOperator{\Ric}{Ric}

\newcommand{\D}{\nabla}
\newcommand{\p}{\partial}

\renewcommand{\bar}{\overline}
\renewcommand{\tilde}{\widetilde}
\renewcommand{\epsilon}{\varepsilon}
\renewcommand{\leq}{\leqslant}
\renewcommand{\geq}{\geqslant}

\newcommand{\TODO}[1]{\message{TODO Warning: {#1}}}
\newcommand{\updatetag}[2]{}

\newtheorem{theorem}{Theorem}[section]
\newtheorem{lemma}[theorem]{Lemma}
\newtheorem{cor}[theorem]{Corollary}

\newtheorem{defn}[theorem]{Definition}
\numberwithin{equation}{section}
\theoremstyle{definition}

\theoremstyle{definition}
\newtheorem{remark}[theorem]{Remark}
\theoremstyle{definition}
\newtheorem{example}[theorem]{Example}

\DeclareMathOperator{\isp}{sip}
\DeclareMathOperator{\ip}{ip}

\renewcommand{\H}{\mathcal{H}}

\DeclareMathOperator{\Lip}{Lip}
\DeclareMathOperator{\spt}{spt}

\begin{document}
	
\title{Isoperimetry and the properness of weak inverse mean curvature flow}
\author{Kai Xu}
\date{}
\maketitle

\begin{abstract}
	 We prove a new existence theorem for proper solutions of Huisken and Ilmanen's weak inverse mean curvature flow, assuming certain non-degeneracy conditions on the isoperimetric profile. In particular, no curvature assumption is imposed in our existence theorem.
\end{abstract}
	
	
%
	
\section{Introduction}


\TODO{the introduction in the submitted version is different from here. Remember to revise accordingly.}

The inverse mean curvature flow (IMCF) is an extrinsic parabolic flow in which a family of hypersurfaces $\Sigma_t$ evolve at the speed of the inverse of mean curvature:
\begin{equation}\label{eq-intro:smooth_flow}
	\frac{\p\Sigma_t}{\p t}=\frac1H\nu.
\end{equation}
Here $H,\nu$ denote the mean curvature and outer normal vector of $\Sigma_t$. It is often convenient to switch to the level set description: if one defines a function $u$ such that $\Sigma_t=\{u=t\}$, then \eqref{eq-intro:smooth_flow} is equivalent to $u$ being a solution of the following degenerate elliptic equation:
\begin{equation}\label{eq-intro:level_set}
	\div\Big(\frac{\D u}{|\D u|}\Big)=|\D u|.
\end{equation}

Starting from a star-shaped strictly mean convex hypersurface in $\RR^n$, Gerhardt \cite{Gerhardt_1991} and Urbas \cite{Urbas_1990} showed that the solution to \eqref{eq-intro:smooth_flow} exists for all $t$ and converges to round spheres as $t\to\infty$. On general Riemannian manifolds or with general initial hypersurfaces, however, the smooth IMCF may inevitably form a finite-time singularity and stop evolving (see for example \cite{Harvie_2022}). In the foundational work of Huisken and Ilmanen \cite{Huisken-Ilmanen_2001}, a weak formulation of the IMCF is developed by finding a suitable variational principle for the level set equation \eqref{eq-intro:level_set}. Heuristically speaking, the weak formulation combines the evolution by $1/H$ with additionally a jumping mechanism. Along the weak flow, there is a set of time $t$ when the surface $\Sigma_t=\p\{u<t\}$ immediately ``jumps'' to $\Sigma_t^+=\p\{u\leq t\}$ which is characterized as the outermost minimizing hull of $\Sigma_t$. The level sets $\Sigma_t$ of a weak solution have $C^{1,1}$ regularity \cite{Heidusch_thesis}, and one can view jumping as the way to bypass the formation of singularity in IMCF.

In \cite{Huisken-Ilmanen_2001} the weak IMCF is applied to resolve the Riemannian Penrose inequality for the case of a single horizon (another proof of the general case is due to Bray \cite{Bray_2001}). The key fact proved in \cite{Huisken-Ilmanen_2001} is that the monotonicity of Hawking mass, initially observed by Geroch \cite{Geroch_1973} in the smooth case, continues to hold along the weak flow. Other applications of the weak IMCF include the computation of Yamabe invariants \cite{Bray-Neves_2004}, isoperimetric inequalities related to scalar curvature \cite{Bray-Miao_2008, Brendle-Chodosh_2014, Shi_2016}, and a recent new proof of the Ricci pinching conjecture \cite{Huisken-Koerber_2023}. In a recent work \cite{Xu_2023} of the author on geometric inequalities for 3-manifolds with positive scalar curvature, the weak IMCF is applied on the universal cover of a compact manifold.

In this paper we focus on the existence of proper weak solutions of IMCF. A weak solution $u$ is called \textit{proper} if all the sub-level sets $E_t=\{u<t\}$ are precompact, or equivalently, if $\lim_{d(x,x_0)\to\infty}u(x)=+\infty$ where $x_0$ is a fixed basepoint. Many desirable properties of weak solutions, including the maximum principle and monotonicity of Hawking mass, require a certain properness of the weak solutions in question. The existence of proper weak solutions is a non-trivial question relating to the geometry of the manifold at infinity. For instance, there does not exist proper solutions on a manifold with a cylindrical end. This is because the level set $\Sigma_t$ immediately escapes to infinity once it entirely enters the cylindrical region. With more work, we can show that the existence of proper solutions implies superlinear growth of volume. The latter suggests that the manifold is expanding at infinity in some sense. In \cite{Huisken-Ilmanen_2001} Huisken and Ilmanen proved an existence theorem for proper weak solutions, assuming the presence of a proper weak subsolution. The proof in \cite{Huisken-Ilmanen_2001} involves an elliptic regularization process for the equation \eqref{eq-intro:level_set}, where the weak subsolution serves as a priori barrier. While finding proper weak subsolutions is straightforward for many spaces with precise description at infinity, such as asymptotically flat manifolds, it remains a difficult task in general.

%
%
%

More recently, the $p$-harmonic approximation scheme introduced by Moser \cite{Moser_2007} has provided another fundamental approach to solving the weak IMCF. Given a smooth bounded region $\Omega\subset\RR^n$, viewed as the interior of the initial hypersurface $\Sigma_0$, one considers the unique $p$-harmonic function $u_p$ with $u_p|_{\p\Omega}=1$ and $\lim_{x\to\infty}u_p(x)=0$. Such function is often called the $p$-capacitary potential. In \cite{Moser_2007} it is proved that the functions $(1-p)\log u_p$ converge as $p\to1$ to the (unique) proper weak solution IMCF. The approach of Moser was further pursued by Kotschwar-Ni \cite{Kotschwar-Ni_2009} and Mari-Rigoli-Setti \cite{Mari-Rigoli-Setti_2022} in various settings. In the latter work, the existence theorem is extended to the following generality (see Theorem 1.3, 1.4, 1.7 there): if a non-compact manifold $M$ satisfies a uniform Ricci curvature lower bound $\Ric\geq-\lambda g$ and a Euclidean-type Sobolev inequality $||f||_{n/(n-1)}\leq C||\D f||_1$, $\forall f\in C^1_0(M)$, then with any precompact smooth initial domain there exists a proper weak solution.

\vspace{12pt}

In the present work we prove independently an existence theorem, which does not require any uniform curvature assumption, and links the properness of weak solutions to the expanding of $M$ at infinity in the sense given by an isoperimetric inequality.

\begin{defn}\label{def-intro:ip}
	For a Riemannian manifold $(M,g)$, define its isoperimetric profile by
	\[\ip(v)=\inf\big\{|\p^*E|: E\subset\subset M\text{ has finite perimeter},\ |E|=v\big\}.\]
\end{defn}

Our main theorem states as follows:

\begin{theorem}\label{thm-intro:existence}
	Let $M$ be a complete, connected, non-compact Riemannian manifold with infinite volume, such that the isoperimetric profile satisfies
	\begin{equation}\label{eq-intro:nondeg_ip}
		\liminf_{v\to\infty}\ip(v)=\infty
		\ \ \ \text{and}\ \ \ 
		\int_0^{v_0}\frac{dv}{\ip(v)}<\infty\,\ \text{for some}\ v_0>0.
	\end{equation}
	Then with any smooth initial domain $E_0$ there exists a unique proper weak solution of the inverse mean curvature flow, in the sense of Definition \ref{def-imcf:IVP}.
\end{theorem}

The first condition in \eqref{eq-intro:nondeg_ip} implies superlinear volume growth and excludes cylindrical ends. This is consistent with the observations made above. The second condition is less transparent at first glance. However, with some work one can show that it implies a uniform lower bound on the volume of geodesic balls $B(x,1)$. This rules out ends of finite volume, which is another obstruction to the properness of weak solutions. We refer to Remark \ref{rmk-sip:main_conditions} for more details on these observations.

For a geometrically controlled manifold (for example, a manifold with $\Ric\geq -\lambda g$ and $\inf_{x\in M}|B(x,1)|>0$), one is usually able to prove $\liminf_{v\to0}\ip(v)v^{-(n-1)/n}>0$; see \cite{Coulhon-Saloff-Coste_1993} and \cite[Theorem 1.3]{Antonelli-Pasqualetto-Pozzetta-Semola_2022b}. In the case with nonnegative sectional curvature and uniform lower bound on $|B(x,1)|$, more conclusions can be obtained. Indeed, from \cite{Antonelli-Pasqualetto-Pozzetta-Semola_2022a, Antonelli-Pozzetta_2023} one sees that the first condition in \eqref{eq-intro:nondeg_ip} is equivalent to superlinear volume growth. Condition \eqref{eq-intro:nondeg_ip} in particular covers the case of Euclidean isoperimetric inequality
\begin{equation}\label{eq-intro:euc_isop_ineq}
	\ip(v)\geq cv^{(n-1)/n}\qquad\text{($\forall\,v>0$)},
\end{equation}
therefore extends the result in \cite[Theorem 1.3]{Mari-Rigoli-Setti_2022} as mentioned above. On the other hand, one cannot obtain uniform gradient estimate without curvature assumptions. Another case of application is when $M$ is a normal covering space of a closed manifold, whose deck transformation group has at least quadratic growth. In this case there holds $\ip(v)\geq cv^{(n-1)/n}$ for small $v$ and $\ip(v)\geq cv^{1/2}$ for large $v$, by a theorem in \cite{Coulhon-Saloff-Coste_1993}.

In the proof, a uniform diameter bound for the sub-level sets $E_t=\{u<t\}$ is obtained, which is given explicitly by \eqref{eq-exist:bound_main_thm}. When $M$ only satisfies $\liminf_{v\to\infty}\ip(v)>A$ for some $A>0$, for example when $M$ has a cylindrical end, we are able to show that there is a weak solution that remains proper until $|\Sigma_t|=A$. See Theorem \ref{thm-exist:main_thm} for the full statement.

It is worth mentioning that when $M$ satisfies the Euclidean isoperimetric inequality \eqref{eq-intro:euc_isop_ineq}, our quantitative bound \eqref{eq-exist:bound_main_thm} reads $\diam(E_t)\leq Ce^{nt/(n-1)}$, hence Theorem \ref{thm-intro:existence} provides solutions with the growth $u(x)\geq\frac{n-1}{n}\log|x|-C$. On the other hand, the growth $u(x)\geq(n-1)\log|x|-C$ is obtained in many of the known results \cite{Huisken-Ilmanen_2001, Mari-Rigoli-Setti_2022, Moser_2007}. It is interesting to ask whether the latter optimal growth is can be proved assuming \eqref{eq-intro:euc_isop_ineq} solely.


\vspace{12pt}

Finally, we present some explanations on the proof of Theorem \ref{thm-intro:existence}. First we note that the $p$-harmonic approach in \cite{Kotschwar-Ni_2009, Mari-Rigoli-Setti_2022, Moser_2007} fails to apply here, since there exist manifolds satisfying \eqref{eq-intro:nondeg_ip} but are $p$-parabolic for all $p>1$. For example, consider the manifold $S^{n-1}\times\RR$ with the metric
\[g=dr^2+\log^2(2+r^2)g_{S^{n-1}},\]
which by \cite[Proposition 1.7]{Holopainen_1999} is $p$-parabolic for all $p>1$. In order to prove Theorem \ref{thm-intro:existence}, we investigate the measure-theoretic aspects of the weak IMCF. In particular, we investigate one of the original weak formulations of Huisken and Ilmanen, in which the sub-level sets of a weak solution $u$ locally minimize the functional
\[E\mapsto|\p^*E|-\int_E|\D u|.\]
We use this minimizing property, along with the isoperimetric inequalities \eqref{eq-intro:nondeg_ip}, to obtain a crucial diameter bound for the sub-level sets (see Lemma \ref{lemma-exist:apriori_bound_IMCF}). This bound is viewed as an a propri properness estimate for weak solutions.

Another ingredient in the proof is a conic cutoff method, inspired by similar arguments apprearing in \cite[Theorem 3.1]{Huisken-Ilmanen_2001}. Given an initial value $E_0$, we let $E_0\subset W_1\subset W_2\subset\cdots$ be an exhaustion of smooth precompact domains. For each $k$ we construct another manifold $M_k$ by appropriately attaching a metric cone $\p W_k\times[0,\infty)$ to $W_k$. The original existence theorem of Huisken and Ilmanen guarantees a weak solution $u_k$ on each $M_k$. Then we show that the value of $u_k$ eventually stabilizes, and that the limit function $u=\lim_{k\to\infty}u_k$ is a well-defined proper weak solution on $M$. In this process, we utilize the diameter estimate in the previous step, as well as the maximum principle for weak IMCF which is of semi-local feature (see Theorem \ref{thm-imcf:max_principle}).

The technique in the present paper may also be applied to the weak IMCF with free boundary, which is developed by Marquardt \cite{Marquardt_2017} and Koerber \cite{Koerber_2020}. \\

\textbf{Organization of the paper}. In Section \ref{sec:imcf} we introduce Huisken and Ilmanen's weak formulation of the inverse mean curvature flow, and review some useful properties of weak solutions. Section 3 contains further preparatory results, where we introduce the notion of strong isoperimetric ratio and strictly outer minimizing hull, and prove a theorem on the existence of precompact minimizing hull. Finally, in Section \ref{sec:exist} we prove the main theorem and its quantitatively refined version, Theorem \ref{thm-exist:main_thm}. \\

\textbf{Notations}. In this paper the manifold $M$ is always assumed to be connected and oriented, and with dimension $n\geq2$. The volume of a subset $X\subset M$ is denoted by $|X|$, and the $k$-dimensional Hausdorff measure is denoted by $\H^k(X)$. For a set $E$ with locally finite perimeter, we denote by $\p^*E$ its reduced boundary. Given an open set $U$, $|\p^*E\cap U|$ denotes the total perimeter of $E$ in $U$. In particular, $|\p^*E|$ is the total perimeter of $E$ in $M$. A subset $K\subset M$ is called precompact if its closure is compact. We use $B(x,r)$ to denote geodesic balls. \\

\textbf{Acknowledgements}. The author would like to thank Gioacchino Antonelli, Mattia Fogagnolo and Thomas Koerber for helpful conversations.
	
\section{The weak inverse mean curvature flow}\label{sec:imcf}

In this section we give a brief introduction to the weak formulation of inverse mean curvature flow, as established by Huisken and Ilmanen \cite{Huisken-Ilmanen_2001}. For more comprehensive introductions to this subject, we refer the reader to Chapter 4 of the book by Lee \cite{Lee_rel}, along with the original work \cite{Huisken-Ilmanen_2001}. The starting point of the theory is the following variational characterization of \eqref{eq-intro:level_set}:

\begin{defn}[weak solution, cf. {\cite[p.365]{Huisken-Ilmanen_2001}}]\label{def-imcf:weak_sol_u} {\ }
	
	For a precompact domain $K$ and two locally Lipschiz functions $u, v$, we define the following energy functional:
	\begin{equation}\label{eq-imcf:energy_u}
		J_u^K(v)=\int_K\big(|\D v|+v|\D u|\big).
	\end{equation}
	Given a domain $\Omega\subset M$. A function $u\in\Lip_{loc}(\Omega)$ is called a weak solution (resp. subsolution, supersolution) of IMCF in $\Omega$, if
	\[J_u^K(u)\leq J_u^K(v)\]
	for all $K\subset\subset\Omega$ and all $v\in\Lip_{loc}(\Omega)$ (resp. all $v\leq u$, $v\geq u$) such that $\{u\ne v\}\subset K$.
\end{defn}

By \cite[Lemma 4.30]{Lee_rel}, a smooth function $u$ with non-vanishing gradient solves \eqref{eq-intro:level_set} if and only if $u$ satisfies Definition \ref{def-imcf:weak_sol_u}. Another weak formulation is introduced in \cite{Huisken-Ilmanen_2001}, which reflects the global minimizing property of level sets. Throughout this paper, we fix the following notations for sub-level sets:

\begin{defn}\label{def-imcf:level_sets}
	For a function $u$, define
	\[E_t(u)=\{u<t\},\qquad E_t^+(u)=\{u\leq t\}.\]
	When there is no ambiguity, we make implicit the dependence on $u$ and write $E_t,E_t^+$.
\end{defn}

The second weak formulation is as follows.

\begin{defn}[weak solution II, cf. {\cite[p.366]{Huisken-Ilmanen_2001}}]\label{def-imcf:weak_sol_E} {\ }
	
	For a precompact domain $K$, a locally Lipschitz function $u$ and a set $E$ of locally finite perimeter, define the following functional:
	\begin{equation}\label{eq-imcf:energy_E}
		J_u^K(E)=|\p^*E\cap K|-\int_{E\cap K}|\D u|.
	\end{equation}
	Given a domain $\Omega\subset M$. A function $u\in\Lip_{loc}(\Omega)$ is called a weak solution (resp. subsolution, supersolution) of IMCF in $\Omega$, if
	\[J_u^K(E_t)\leq J_u^K(E)\]
	for all $t\in\RR$, all $K\subset\subset\Omega$ and all set $E$ (resp. all $E\supset E_t$, $E\subset E_t$) with locally finite perimeter, such that $E\Delta E_t\subset K$.
\end{defn}

A direct but important consequence of this definition is that: if $u$ is a weak solution in $M$, then each sub-level set $E_t(u)$ is outward minimizing in $M$ in the sense of Definition \ref{def-imcf:outer_area_minimizer} below, whenever it is precompact. This suggests that jumping occurs when $E_t$ ceases to be outward minimizing.

\begin{defn}[outward area minimizer]\label{def-imcf:outer_area_minimizer} {\ }
	
	A precompact set $E\subset M$ with finite perimeter is said to be outward minimizing (resp. strictly outward minimizing), if $|\p^*E|\leq |\p^*F|$ (resp. $|\p^*E|<|\p^*F|$) for any other precompact set $F\supset E$ with $|F\setminus E|>0$.
\end{defn}

In \cite[Lemma 1.1]{Huisken-Ilmanen_2001} one proves that Definitions \ref{def-imcf:weak_sol_u} and \ref{def-imcf:weak_sol_E} are equivalent. The following properties of weak solutions are straightforward to verify (see \cite{Huisken-Ilmanen_2001}). (1) If $u$ is a weak solution in $\Omega$ then $u$ is a weak solution in any subdomain $\Omega'\subset\Omega$. (2) $u$ is a weak solution if and only if it is simultaneously a weak subsolution and supersolution. (3) If $u$ is a weak subsolution (resp. supersolution) in $\Omega$, then $\min(u,t)$ is also a weak subsolution (resp. supersolution) in $\Omega$, for any $t\in\RR$. (4) The constant function is trivially a weak solution, and there is no non-constant weak solution on a compact manifold. Furthermore, Definition \ref{def-imcf:weak_sol_E} implies $J_u^K(E_s)=J_u^K(E_t)$ for all $s<t$, whenever $E_t\setminus E_s\subset K\subset\subset\Omega$. By the coarea formula, this implies the exponential growth of area
\[|\p^*E_t|=e^{t-s}|\p^*E_s|.\]

Next we define the initial value problem for the weak IMCF.

\begin{defn}[initial value problem, cf. {\cite[p.367-368]{Huisken-Ilmanen_2001}}]\label{def-imcf:IVP} {\ }
	
	Let $E_0\subset\subset \Omega$ be a precompact $C^{1,1}$ domain. We say that $u$ is a weak solution of IMCF on $\Omega$ with initial condition $E_0$, if
	
	(1) $u\in\Lip_{loc}(\Omega)$ and $E_0=\{u<0\}$,
	
	(2) $u$ is a weak solution to the IMCF in $\Omega\setminus\bar{E_0}$ in the sense of Definition \ref{def-imcf:weak_sol_u} or \ref{def-imcf:weak_sol_E}.
	
	A weak solution $u$ (of the initial value problem) is called \textit{proper}, if $E_t(u)$ is precompact in $M$ for all $t\geq0$, equivalently, if $\lim_{x\to\infty}u(x)=+\infty$ uniformly.
\end{defn}

Note that this definition does not concern the value of $u$ or the geometry in the interior of $E_0$. The following example describes specifically the case of spherical symmetry.

\begin{example}\label{ex-imcf:sph_sym}
	Consider the metric of the form $g=dr^2+f(r)^2g_{S^{n-1}}$, $r\in[0,\infty)$, where $f$ is smooth and positive for $r>0$. For a geodesic ball $E_0=\{r<r_0\}$, one set
	\begin{equation}\label{eq-imcf:rot_sym}
		u(r)=\max\Big\{0,(n-1)\inf_{s\geq r}\log\frac{f(s)}{f(r_0)}\Big\}\qquad\text{for}\ \ r\geq r_0,
	\end{equation}
	and one extends $u$ arbitrarily so that $u<0$ in $E_0$. The resulting function is a weak solution of IMCF with initial condition $E_0$. When $f$ is strictly increasing in $[r_0,\infty)$, the weak solution coincides with the classical solution. When $f$ is not monotone, jumping occurs along the flow. In particular, if $E_0$ is not strictly outward minimizing, (i.e. if there exists $r>r_0$ with $f(r)\leq f(r_0)$), a jumping occurs at $t=0$. An extreme case is $\liminf_{r\to\infty}f(r)\leq f(r_0)$, for which there is no weak solutions other than $u\equiv0$. The maximal time of existence is $T=(n-1)\log\big(\liminf_{r\to\infty}f(r)/f(r_0)\big)$, in the sense that $E_t$ is precompact only when $t<T$. In particular, the solution $u$ is proper if and only if $\liminf_{r\to\infty}f(r)=\infty$. These are particular cases of Theorem \ref{thm-exist:main_thm}.
\end{example}

Generally, at each jumping time the set $E_t^+$ is characterized as the \textit{strictly outward minimizing hull} of $E_t$, which we refer to subsection \ref{subsec:hull} for more details.

The following existence theorem is proved in \cite[Theorem 3.1]{Huisken-Ilmanen_2001}, using the method of elliptic regularization. The only requirement here is the existence of a proper subsolution outside of a compact set. When $M$ is a $C^0$-asymptotically flat manifold, the function $w=(n-1-\epsilon)\log|x|$ serves as a weak subsolution \cite{Benatti-Fogagnolo-Mazzieri_2022}. When $M$ is a Cartan-Hadamard manifold (i.e. a simply-connected manifold with nonpositive sectional curvature), the function $w=(n-1)\log d(x,x_0)$ is a classical subsolution in the sense that $\div\big(\frac{\D w}{|\D w|}\big)\geq|\D w|$, by Hessian comparison.

\begin{theorem}[barrier implies existence]\label{thm-imcf:subsol_implies_exist} {\ }
	
	Let $(M,g)$ be a complete connected Riemannian $n$-manifold without boundary. If there exists a proper, locally Lipschitz, weak subsolution of IMCF with a precompact initial condition, then for any nonempty smooth precompact domain $E_0\subset M$ there exists a proper locally Lipschitz weak solution to IMCF on $M$ with initial condition $E_0$. Such solution is unique among proper solutions on $M\setminus E_0$.
	
	Furthermore, the following local gradient estimate holds for the solution obtained:
	\begin{equation}\label{eq-imcf:gradient_estimate}
		|\D u(x)|\leq \sup_{\p E_0\cap B_r(x)}H_++C(n)r^{-1},\qquad \forall a.e.\ x\in M\setminus E_0,
	\end{equation}
	for each $0<r<\sigma(x)$, where $H_+$ is the positive part of the mean curvature of $\p E_0$, and $\sigma(x)$ is a constant depending on the geometry near $x$, see \cite[Definition 3.3]{Huisken-Ilmanen_2001}.
\end{theorem}

Finally, we summarize several properties of weak solutions.

\begin{lemma}[properties of level sets]\label{lemma-imcf:Et_properties} {\ }
	
	Let $u$ be a weak solution to the IMCF on $M$ with smooth initial condition $E_0$. Then:
	
	(1) Each $\p E_t$ and $\p E_t^+$ is a $C^{1,\alpha}$ ($\alpha<\frac12$) hypersurface except for a singular set of Hausdorff dimension at most $n-8$. Each $E_t$ ($t>0$) is outward area-minimizing, and each $E_t^+$ ($t\geq0$) is strictly outward minimizing, provided $E_t$ resp. $E_t^+$ is precompact.
	
	(2) For all $t\geq0$, $E_t^+$ is the strict outward minimizing hull of $E_t$ in the sense of Definition \ref{def-hull:least_volume_env}, provided that $E_t^+$ is precompact.
	
	(3) For all $t>0$, we have $\p E_s\to \p E_t$ when $s\nearrow t$ and $\p E_s\to \p E_t^+$ when $s\searrow t$: the convergence holds in the local $C^{1,\beta}$ topology ($\beta<\alpha$) when $n\leq7$, and in local Hausdorff topology in all dimensions.
	
	(4) For any $t>0$ we have $|\p E_t|=|\p E_t^+|$, provided the latter is precompact. For $t=0$ we have $|\p E_0|\geq|\p E_0^+|$ provided the latter is precompact, and equality holds if $E_0$ is outward minimizing.
	
	(5) The area of level sets is increasing exponentially: $|\p E_t|=e^t|\p E_0^+|$ whenever the former is precompact, for all $t>0$.
\end{lemma}

\begin{proof}
	Items (1)\,$\sim$\,(4) follows from Definition \ref{def-imcf:weak_sol_E}, \cite[p.368-369]{Huisken-Ilmanen_2001}, \cite[Property 1.4]{Huisken-Ilmanen_2001}, and the fundamental results in geometric mesure theory (see \cite[Chapter 21]{Maggi}). Part (5) follows from $J_u(E_{t_1})=J_u(E_{t_2})$, $\forall\,0<t_1<t_2$, see also \cite[Lemma 1.6]{Huisken-Ilmanen_2001}.
\end{proof}

\begin{theorem}[maximum principle, {\cite[Theorem 2.2]{Huisken-Ilmanen_2001}}]\label{thm-imcf:max_principle} {\ }
	
	(1) If $u,v$ are two weak solutions to the IMCF in a domain $\Omega$, such that $\{u\ne v\}\subset\subset\Omega$, then $u=v$ in $\Omega$.
	
	(2) From a given precompact initial condition $E_0$ there exists at most one proper weak solution. More generally, if $u_1,u_2$ are two weak solutions in $\Omega$ with the same initial condition $E_0\subset\subset\Omega$, then $E_t(u_1)=E_t(u_2)$ whenever $E_t(u_1)$, $E_t(u_2)$ are precompact in $\Omega$.
\end{theorem}


\section{Isoperimetric profile and minimizing hull}\label{sec:sip}

We include in this section some discussions relevant to the isoperimetric profile, primarily in preparation for the proof of Theorem \ref{thm-intro:existence}. Some results are also of their own interests.

\begin{remark}
	We will be frequently using the following fact: if $E,F\subset M$ are precompact sets with bounded perimeter, then
	\begin{equation}\label{eq-sip:cap_cup}
		|\p^*(E\cap F)|+|\p^*(E\cup F)|\leq|\p^*E|+|\p^*F|.
	\end{equation}
	See \cite[Lemma 12.22]{Maggi} for a proof.
	The following decomposition lemma is also useful, see \cite[Theorem 16.3]{Maggi} for its proof:
	
	\begin{lemma}\label{lemma-sip:intersect_decomp}
		Let $E,F\subset M$ be two precompact sets of finite perimeter. Assuming that $\H^{n-1}(\p^*E\cap\p^*F)=0$, we have
		\begin{equation}\label{eq-sip:intersect_decomp}
			|\p^*(E\cap F)|=\H^{n-1}\big(\p^*E\cap F^{(1)}\big)+\H^{n-1}\big(\p^*F\cap E^{(1)}\big),
		\end{equation}
		where $E^{(1)}$ is the measure-theoretic interior of $E$, defined as
		\[E^{(1)}=\Big\{x\in M:\lim_{r\to0}\frac{|E\cap B(x,r)|}{|B(x,r)|}=1\Big\}.\]
	\end{lemma}
	
	Note that $|E\Delta E^{(1)}|=0$ from Lebesgue's differentiation lemma. It is easy to see that \eqref{eq-sip:intersect_decomp} does not hold if $F^{(1)}$ (or $E^{(1)}$) is replaced by $F$ (or $E$). When $E$ is a domain with $C^1$ boundary, we have $E^{(1)}=E$ pointwise.
	
	By Federer's structure theorem \cite[Theorem 16.2]{Maggi}, Lemma \ref{lemma-sip:intersect_decomp} implies
	\begin{equation}\label{eq-sip:intersection_peri}
		|\p^*E|-|\p^*(E\cap F)|=\H^{n-1}(\p^*E\cap(M\setminus F)^{(1)})-\H^{n-1}(\p^*F\cap E^{(1)}),
	\end{equation}
	whenever $\H^{n-1}(\p^*E\cap\p^*F)=0$.
\end{remark}

\begin{remark}\label{rmk-sip:geo_sph}
	Let $M$ be a complete manifold. Consider the family of geodesic balls $B_r=B(x_0,r)$. By the main theorem of \cite{Itoh-Tanaka_1998}, for almost every $r\in\RR$ the geodesic sphere $S_r=\p B_r$ is a smooth hypersurface except at a singular set of zero $(n-1)$-Hausdorff measure. For these $r$, both $B_r$ and $M\setminus\bar{B_r}$ coincide with their measure-theoretic interior up to zero $(n-1)$-Hausdorff measure.
\end{remark}

\begin{remark}\label{rmk-sip:level_set}
	If $u$ is a weak solution of IMCF in a domain $\Omega$, then we have $E_t=E_t^{(1)}$ for each $t$. Since $E_t$ is open, this claim will follow from $\p E_t\cap E_t^{(1)}=\emptyset$. Fix $x\in\p E_t$. Since $u$ is continuous, we may find a sequence $t_i\nearrow t$ and $x_i\to x$, such that $x_i\in\p^*E_{t_i}$. By Definition \ref{def-imcf:weak_sol_E} and \cite[Theorem 21.14]{Maggi}, we obtain $x\in\spt|D\chi_{E_t}|$. Then by the minimizing effect of $E_t$, we obtain $x\notin E_t^{(1)}$, see \cite[Theorem 21.11]{Maggi}. This proves our claim. Similarly, we can show that $\Omega\setminus E_t^+=(\Omega\setminus E_t^+)^{(1)}$. Finally, note that these facts imply $\p E_t=\spt|D\chi_{E_t}|$ and $\p E_t^+=\spt|D\chi_{E_t^+}|$.
\end{remark}

\subsection{The strong isoperimetric profile}

\begin{defn}\label{def-sip:sip}
	Given a complete manifold $M$, define the strong isoperimetric profile
	\[\isp(v)=\inf\big\{|\p^*E|: E\subset\subset M, |E|\geq v\big\},\]
	for $0<v<|M|$. The inverse isoperimetric profile, denoted by $\isp^{-1}$ for convenience, is defined as
	\[\isp^{-1}(a)=\sup\big\{|E|: E\subset\subset M, |\p^*E|\leq a\big\}.\]
\end{defn}

Clearly $\isp(v)=\inf_{v'\geq v}\ip(v')$; where we recall Definition \ref{def-intro:ip} for the function $\ip(v)$. The strong isoperimetric profile has less pathologic behavior to be concerned. In particular, $\isp(v)$ is always continuous when $M$ is connected, while the continuity of $\ip(v)$ is a delicate problem. See \cite{Flores-Nardulli_2019, Nardulli-Pansu_2018} and more recently \cite[Corollary 4.16]{Antonelli-Pasqualetto-Pozzetta-Semola_2022a} which address the latter problem. On the other hand, the strong isoperimetric profile is a non-trivial quantity only for manifolds that are non-degenerate at infinity. For example, a manifold with finite volume has $\isp(v)\equiv0$. It is easy to examine that both $\isp(v)$ and $\isp^{-1}(a)$ are non-decreasing, and $\isp^{-1}$ is defined with finite value on the interval $\big(0,\lim_{v\to\infty}\isp(v)\big)=\big(0,\liminf_{v\to\infty}\ip(v)\big)$.


\begin{lemma}\label{lemma-sip:cont}
	Suppose $M$ is complete, connected, and has infinite volume. Then $\isp(v)$ is continuous on $(0,\infty)$.
\end{lemma}
\begin{proof}
	Since $\isp(v)$ is non-decreasing, its right continuity at a value $v$ is equivalent to
	\begin{equation}\label{eq-sip:right_cont}
		\inf_{|E|\geq v}|\p^*E|\geq\inf_{|E|>v}|\p^*E|.
	\end{equation}
	Given any precompact set $E$ with $|E|\geq v$ and $|\p^*E|\leq\isp(v)+\epsilon$, we choose any geodesic ball $B$ of perimeter at most $\epsilon$, such that $|B\setminus E|>0$. Then note that $|E\cup B|>v$ and $|\p^*(E\cup B)|\leq\isp(v)+2\epsilon$. Letting $\epsilon\to0$ this proves \eqref{eq-sip:right_cont}. It remains to prove left continuity at any fixed $v>0$. Denote
	\[a:=\lim_{v'\to v^-}\isp(v')=\sup_{v'<v}\isp(v').\]
	Fix a basepoint $x_0\in M$. Let $R$ be sufficiently large such that $|B(x_0,R)|>2\cdot 10^nv$. For any $\epsilon>0$, there exists $\delta>0$ such that: for all $x\in B(x_0,R)$, the unique geodesic ball centered at $x$ and with volume $\delta$ has perimeter $<\epsilon$. Denote this geodesic ball by $B(x,r_x)$. By further decreasing $\delta$, we may assume that $|B(x,5r_x)|\leq 10^n|B(x,r_x)|$ for all $x\in B(x_0,R)$. By the definition of $\isp(v)$, for any $v-\delta/2<v'<v$ there exists a precompact set $E$ with $|E|\geq v'$ and $|\p^*E|\leq\isp(v')+\epsilon\leq a+\epsilon$. If $|E|\geq v$, then this already implies $\isp(v)\leq a+\epsilon$. Now assume $|E|\leq v$. We claim that there is a point $x\in B(x_0,R)$ such that $|E\cap B(x,r_x)|\leq\frac12|B(x,r_x)|$. Otherwise, we would have $|E\cap B(x,r_x)|>\frac12|B(x,r_x)|$ for all $x\in B(x_0,R)$. By Vitali's covering lemma, there is a disjoint countable collection of balls $B(x_i,r_{x_i})$ such that $B(x_i,5r_{x_i})$ covers $B(x_0,R)$. Thus
	\[|E|\geq\sum|E\cap B(x_i,r_{x_i})|\geq\frac12\cdot 10^{-n}\sum|B(x_i,5r_{x_i})|\geq\frac12\cdot 10^{-n}|B(x_0,R)|>v,\]
	which is a contradiction. After finding the ball $B(x,r_x)$, consider the new set $E'=E\cup B(x,r_x)$. We have $|E'|\geq|E|+\frac12|B(x,r_x)|>(v-\delta/2)+\delta/2\geq v$ and $|\p^*E'|\leq|\p^*E|+|\p B(x,r_x)|\leq a+2\epsilon$, hence $\isp(v)\leq a+2\epsilon$. Thus either $|E|\geq v$ or $|E|\leq v$ we have obtained $\isp(v)\leq a+2\epsilon$. Taking $\epsilon\to0$ this implies $\isp(v)\leq a$, which proves the left continuity.
\end{proof}

\begin{lemma}\label{lemma-sip:isp_nonzero}
	Suppose $M$ is complete, connected, and has infinite volume. Then either $\isp\equiv 0$ or $\isp(v)>0$ for all $v>0$.
\end{lemma}
\begin{proof}
	Suppose $\isp(v)$ is not identically zero but is zero somewhere. By monotonicity and continuity, there exists a value $v_0>0$ such that $\isp(v_0)=0$ but $\isp(v)>0$ for all $v>v_0$. By definition, there is a sequence of bounded sets $E_k$ such that $|E_k|\geq v_0$ and $|\p^*E_k|\leq 1/k$. By the maximality of $v_0$, we can assume $|E_k|\leq\frac32v_0$ for all $k$. For a fixed $k\in\NN$, $E_k$ is contained in some bounded connected smooth domain $\Omega$. Enlarging $\Omega$ if necessary, we may assume that $|\Omega|\geq3v_0$. Note that this selection of $\Omega$ requires the connectedness and infinite volume of $M$. By relative isoperimetric inequality, for all $l>k$ we have $|E_l\cap\Omega|\leq C(\Omega)|\p^*E_l\cap\Omega|^{n/(n-1)}\leq C(\Omega)l^{-n/(n-1)}$, therefore $|E_l\cap\Omega|\leq v_0/2$ for sufficiently large $l$. Now the set $E_k\cup E_l$ has volume $\geq3v_0/2$ and perimeter $\leq 2/k$, hence $\isp(3v_0/2)\leq2/k$. Letting $k\to\infty$ this yields $\isp(3v_0/2)=0$, which contradicts the maximality of $v_0$.
\end{proof}

Finally, we prove several claims made in the introduction.

\begin{remark}[on the main condition \eqref{eq-intro:nondeg_ip}]\label{rmk-sip:main_conditions}
	Let $M$ be complete with infinite volume.
	
	First, we note that the condition $\liminf_{v\to\infty}\ip(v)=\infty$ implies superlinear volume growth. Indeed, if $|B(x_0,r)|\leq Cr$ for some constant $C$ and for all $r>1$, then applying the coarea formula, we can find a sequence of radii $r_i\in[2^{i-1},2^i]$ such that $|\p^* B(x_0,r_i)|\leq 2C$. This directly implies $\liminf_{v\to\infty}\ip(v)\leq2C$.
	
	Next, we note that \eqref{eq-intro:nondeg_ip} implies a uniform lower bound for $|B(x,1)|$, for all $x$. To see this, we fix $x\in M$ and denote $V(r)=|B(x,r)|$. Then the coarea formula and isoperimetric inequality implies $dV/dr\geq\ip(V)$. By Lemma \ref{lemma-sip:isp_nonzero} and \eqref{eq-intro:nondeg_ip}, the quantitiy $\int_0^V dv/\ip(v)$ is defined and finite for all $v$. Hence the ODE inequality can be integrated to give $\int_0^{V(1)}dv/\ip(v)\geq1$, which is the desired uniform lower bound on $|B(x,1)|$.
\end{remark}

\subsection{Strictly outward minimizing hull}\label{subsec:hull}

At each jumping time of the weak IMCF, the set $E_t^+$ is characterized as the strictly outward area-minimizing hull of $E_t$. In \cite{Huisken-Ilmanen_2001} this is described as the smallest strictly outward minimizing set that contains $E_t(u)$ (see Definition \ref{def-hull:least_volume_env}). The theory of minimizing hulls is developed systematically by Fogagnolo and Mazzieri \cite{Fogagnolo-Mazzieri_2022}, where two equivalent definitions are proposed. We briefly introduce the two formulations following \cite{Fogagnolo-Mazzieri_2022}, and prove an equivalence result that improves \cite[Theorem 2.18]{Fogagnolo-Mazzieri_2022}. For the statements in this subsection, the manifold $M$ need not be complete.

\begin{defn}[least area problem; {\cite[Definition 2.6]{Fogagnolo-Mazzieri_2022}}]\label{def-hull:max_vol_sol} {\ }
	
	Given a precompact domain $\Omega\subset\subset M$ with finite perimeter. We say that a precompact set $E\supset\Omega$ solves the least area problem outside $\Omega$, if
	\[|\p^*E|=\inf\big\{|\p^*F|:\Omega\subset F\subset\subset M\big\}.\]
	We say that a precompact set $E$ is a maximal volume solution to the least area problem, if
	\[|E|=\sup\big\{|F|: F\text{ solves the least area problem outside }\Omega\big\}.\]
\end{defn}

The following facts are straightforward to verify, by repeatedly using \eqref{eq-sip:cap_cup}.

\begin{lemma}\label{lemma-hull:max_vol_sol}
	If $E_1,E_2$ both solve the least area problem outside $\Omega$, so does $E_1\cup E_2$. Therefore, the maximal volume solution (to the least area problem) is unique up to measure zero, if it exists. The maximal volume solution is always strictly outward minimizing. If $E$ is the maximal volume solution, and $E'$ is another solution to the least area problem, then $E'\subset E$ up to measure zero.
\end{lemma}

The other definition of minimizing hull is characterized from outside. For the notion of strictly outward minimizing sets, see Definition \ref{def-imcf:outer_area_minimizer}.

\begin{defn}[least volume strictly minimizing envelope; {\cite[Definition 2.12]{Fogagnolo-Mazzieri_2022}}]\label{def-hull:least_volume_env} {\ }
	
	Given a precompact domain $\Omega\subset\subset M$ with finite perimeter, denote
	\[\mathcal{F}(\Omega)=\big\{F: \Omega\subset F\subset\subset M\text{ and }F\text{ is strictly outward minimizing}\big\}.\]
	A precompact set $E\supset\Omega$ is called a least volume strictly minimizing envelope of $\Omega$, if $\mathcal{F}(\Omega)$ is non-empty and
	\[|E|=\inf_{F\subset\mathcal{F}(\Omega)}|F|.\]
\end{defn}

It is not hard to verify the following (see the proof of \cite[Theorem 2.15]{Fogagnolo-Mazzieri_2022}):

\begin{lemma}
	If two sets $E_1,E_2$ are both strictly outward minimizing, then so does $E_1\cap E_2$. Therefore, the least volume strictly minimizing envelope is unique up to measure zero, if exists.
\end{lemma}

In \cite[Theorem 2.16]{Fogagnolo-Mazzieri_2022} the following three statements are proved to be equivalent: (1) every precompact $\Omega$ admits a precompact maximal volume solution as in Definition \ref{def-hull:max_vol_sol}, (2) every precompact $\Omega$ admits a precompact least volume minimizing envelope as in Definition \ref{def-hull:least_volume_env}, and (3) there is an exhaustion of $M$ by precompact strictly outward minimizing hulls. Here we observe that the equivalence between (1) and (2) holds at the level of each single $\Omega$.

\begin{theorem}[equivalence of the two formulations]\label{def-hull:minimizing_hull} {\ }
	
	Given a precompact domain $\Omega\subset M$ with finite perimeter. The following statements are equivalent:
	
	(1) There exists a precompact maximal volume solution $E_1$ to the least area problem outside $\Omega$, in the sense of Definition \ref{def-hull:max_vol_sol}.
	
	(2) There exists a precompact least volume strictly minimizing envelope $E_2$ of $\Omega$, in the sense of Definition \ref{def-hull:least_volume_env}.
	
	Moreover, we have $E_1=E_2$ up to measure zero, if either of them exists. We call $E_1$ the strictly outward minimizing hull (or ``minimizing hull'' for brevity) of $\Omega$.
\end{theorem}

\begin{proof}
	First suppose $E_1$ exists as in (1). To prove (2) it is sufficient to show that $E_1$ is a least volume strictly minimizing envelope of $\Omega$. It is not hard to see from the definition that $E_1$ is strictly outward minimizing. Suppose $F\supset\Omega$ is another precompact strictly outward minimizing set. Since $|\p^*(F\cap E_1)|\geq|\p^*E_1|$ by the minimizing property of $E_1$, we have $|\p^*(F\cup E_1)|\leq|\p^*F|$ by \eqref{eq-sip:cap_cup}, which implies $E_1\subset F$ up to measure zero. This proves (2).
	
	Now suppose that $E_2$ exists as in (2), and we prove that $E_2$ is a maximal volume solution to the least area problem outside $\Omega$. We first show that $E_2$ solves the least area problem. Suppose this is not true, so that there is another precompact set $F\supset\Omega$ with $|\p^*F|<|\p^*E_2|$. Since $|\p^*(F\cup E_2)|\geq|\p^*E_2|$ by the minimizing property of $E_2$, we have $|\p^*(F\cap E_2)|\leq|\p^*F|$. Consider the following least area problem with obstacles
	\[A=\inf\big\{|\p^*G|:\Omega\subset G\subset E_2\big\},\]
	and the volume maximization problem among least area sets
	\[V=\sup\big\{|G|:\Omega\subset G\subset E_2,|\p^*G|=A\big\}.\]
	By the classical compactness theorem \cite[Theorem 12.26]{Maggi}, there exists a maximal volume solution $G_0$ with $|G_0|=V$, $|\p^*G_0|=A$ and $\Omega\subset G_0\subset E_2$. We claim that $G_0$ is strictly outward minimizing in $M$. Suppose $H\supset G_0$ is precompact in $M$. By the minimizing property of $E_2$, we have
	\begin{equation}\label{eq-hull:aux1}
		|\p^*(H\cup E_2)|\geq|\p^*E_2|\ \Rightarrow\ |\p^*(H\cap E_2)|\leq|\p^*H|.
	\end{equation}
	By the minimizing property of $G_0$, we have
	\begin{equation}\label{eq-hull:aux2}
		|\p^*G_0|\leq|\p^*(H\cap E_2)|.
	\end{equation}
	These imply that $|\p^*G_0|\leq|\p^*H|$. Moreover, if equality holds then \eqref{eq-hull:aux1} \eqref{eq-hull:aux2} also attains equality. From the strict minimizing property of $G_0$ (in $E_2$) and $E_2$ (in $M$), we have $H=G_0$ up to measure zero. This proves strict outward minimizing of $G_0$. We have now found a strictly outward minimizing set $G_0$ inside $E_2$, with $|\p^*G_0|\leq|\p^*(F\cap E_2)|\leq|\p^*F|<|\p^*E_2|$. Therefore $|E_2\setminus G_0|>0$ strictly, which contradicts the volume minimizing property of $E_2$.
	
	We have proved that $E_2$ solves the least area problem. Suppose there is another precompact $F\supset\Omega$ with $|\p^*F|=|\p^*E_2|$. Since $|\p^*(F\cap E_2)|\geq|\p^*F|$, we conclude $|\p^*(F\cup E_2)|\leq|\p^*E_2|$. This implies $F\subset E_2$ by the minimizing property of $E_2$, hence $E_2$ is the maximal volume solution. The proof of the theorem is complete.
\end{proof}


\subsection{Existence of precompact minimizing hull}

In general, there may not exist a precompact minimizing hull for a given $\Omega\subset M$. This happens, for example, when $M$ has a cylindrical end or an end of finite volume. In \cite{Fogagnolo-Mazzieri_2022} an existence theorem is proved under two assumptions: (1) $M$ satisfies an Euclidean-like isoperimetric inequality, or (2) $M$ has nonnegative Ricci curvature and superlinear volume growth at infinity. Case (1) is derived from an a priori bound, while case (2) is proved by utilizing an existence theorem for the weak IMCF in \cite{Mari-Rigoli-Setti_2022}. Here we prove an extension of case (1) by weakening the isoperimetric condition of $M$ to \eqref{eq-sip:nondeg_ip}. The proof is in similar spirit with \cite[Proposition 3.1]{Fogagnolo-Mazzieri_2022}, and the crucial part is the following a priori bound for area minimizers:

\begin{lemma}\label{lemma-sip:apriori_barrier}
	Let $A>0$ be a constant. Let $M$ be a complete connected manifold with infinite volume. Fix a basepoint $x_0\in M$. Assume that $M$ satisfies
	\begin{equation}\label{eq-sip:nondeg_ip}
		\liminf_{v\to\infty}\ip(v)>A
		\ \ \ \text{and}\ \ \ 
		\int_0^{v_0}\frac{dv}{\ip(v)}<\infty\,\ \text{for some}\ v_0>0.
	\end{equation}
	For any radius $r>0$ and any precompact set $E$ with $|\p^*E|\leq A$, set
	\begin{equation}\label{eq-sip:choice_of_R}
		R=r+1+2\int_0^{\isp^{-1}(|\p^*E|)}\frac{dv}{\ip(v)}.
	\end{equation}
	Then $R$ is finite. Furthermore, the following holds. If $|\p^*E|\leq|\p^*F|$ for all sets $F$ with $E\cap B(x_0,r)\subset F\subset E$, then $E\subset B(x_0,R)$ up to measure zero.
\end{lemma}
\begin{proof}
	The condition \eqref{eq-sip:nondeg_ip} implies $\isp^{-1}(|\p^*E|)\leq\isp^{-1}(A)<\infty$. Lemma \ref{lemma-sip:isp_nonzero} implies $\ip(v)\geq\isp(v_0)>0$ for all $v\geq v_0$, thus $R<\infty$. Suppose on the contrary that $|E\setminus B(x_0,R)|>0$. After a change of measure zero, we may assume $E=E^{(1)}$ (see Lemma \ref{lemma-sip:intersect_decomp}). For almost every $\rho\in[r,R]$ we have $\H^{n-1}(\p^*E\cap\p B(x_0,\rho))=0$ and $\p B(x_0,\rho)$ is smooth up to zero $(n-1)$-Hausdorff measure. For these $\rho$ the following are defined:
	\[V(\rho)=\big|E\setminus B(x_0,\rho)\big|,\quad 
	S(\rho)=\H^{n-1}\big(E\cap\p B(x_0,\rho)\big),\quad
	A(\rho)=\H^{n-1}\big(\p^*E\setminus\bar{B(x_0,\rho)}\big).\]
	By coarea formula, $V(\rho)$ is absolutely continuous and satisfy
	\begin{equation}\label{eq-sip:stable1}
		V'(\rho)=-S(\rho)
	\end{equation}
	for a.e. $\rho$. By \eqref{eq-sip:intersection_peri} and the minimizing property of $E$, we have
	\begin{equation}\label{eq-sip:stable2}
		|\p^*E|\leq|\p^*(E\cap B(x_0,\rho))|
		\ \Rightarrow\ 
		A(\rho)\leq S(\rho).
	\end{equation}
	By Lemma \ref{lemma-sip:intersect_decomp} we have
	\begin{equation}\label{eq-sip:stable3}
		S(\rho)+A(\rho)=\big|\p^*(E\setminus\bar{B(x_0,\rho)})\big|\geq\ip(V(\rho)).
	\end{equation}
	Combining \eqref{eq-sip:stable1} \eqref{eq-sip:stable2} \eqref{eq-sip:stable3} we have
	\[V'(\rho)\leq-\frac12\ip(V(\rho))\quad\text{for a.e. }\rho.\]
	Since $V(R)>0$, we can integrate this inequality to obtain
	\[2\int_{V(R)}^{V(r)}\frac{dv}{\ip(v)}\geq R-r.\]
	On the other hand, we have $V(r)\leq|E|\leq\isp^{-1}(|\p^*E|)$. This contradicts \eqref{eq-sip:choice_of_R}.
\end{proof}

\begin{cor}[existence of precompact minimizing hull]\label{cor-sip:minimizing_hull} {\ }
	
	Assume $M$ satisfies the conditions in Lemma \ref{lemma-sip:apriori_barrier}. Then any precompact domain $\Omega$ with $|\p^*\Omega|\leq A$ admits a precompact strictly outward minimizing hull $E$ in the sense of Definition \ref{def-hull:minimizing_hull}. Moreover, if $\Omega\subset B(x_0,r)$, then $E\subset B(x_0,R)$, where $R$ is given by \eqref{eq-sip:choice_of_R}.
\end{cor}

\begin{proof}
	Suppose $\Omega\subset B(x_0,r)$. Denote
	\begin{equation}\label{eq-sip:aux1}
		A'=\inf\big\{|\p^*E|:\Omega\subset E\subset B(x_0,R+1)\}.
	\end{equation}
	By compactness theorem and lower semi-continuity of perimeter, the least perimeter problem for $A'$ always have a solution: there is a domain $E$ with $\Omega\subset E\subset B(x_0,R+1)$, such that $|\p^*E|=A'$. Then we consider the maximal volume solution: set
	\[V=\sup\big\{|F|:\Omega\subset F\subset B(x_0,R+1), \,|\p^*F|=A'\big\}.\]
	By compactness theorem and lower semi-continuity again, there is a set $F_0$ with $E\subset F_0\subset B(x_0,R+1)$, such that $|\p^*F_0|=A'$ and $|F_0|=V$. Note that $F_0$ satisfies the condition in Lemma \ref{lemma-sip:apriori_barrier}, since $F_0$ solves the area minimization problem \eqref{eq-sip:aux1}. From this we conclude that $F_0\subset B(x_0,R)$. We claim that $F_0$ is the minimizing hull of $E$ in $M$, in the sense of Definition \ref{def-hull:max_vol_sol}. Suppose this is not true, so there is another precompact set $F_1\supset\Omega$, with either $|\p^*F_1|<|\p^*F_0|$ or $|F_1|>|F_0|$. There is some radius $R_1$ such that $F_1\subset B(x_0,R_1)$. Repeating the argument above, there exists a maximal volume solution $F_2$ to the least area problem with inner obstacle $\Omega$ and outer obstacle $B(x_0,R_1+1)$. As $F_1$ is a valid competitor in this minimization problem, we obtain that either $|\p^*F_2|<|\p^*F_1|<|\p^*F_0|$ or $|F_2|>|F_1|>|F_0|$. On the other hand, Lemma \ref{lemma-sip:apriori_barrier} implies that $F_2\subset B(x_0,R)$, so we obtain a contradiction with the minimizing property of $F_0$.
\end{proof}


\section{The main existence theorem}\label{sec:exist}

We prove the following quantitative refinement of Theorem \ref{thm-intro:existence}, which will imply Theorem \ref{thm-intro:existence} by a limiting argument.

\begin{theorem}\label{thm-exist:main_thm}
	Given a constant $A>0$. Let $M$ be a connected, complete, non-compact $n$-manifold ($3\leq n\leq 7$) with infinite volume, such that
	\begin{equation}\label{eq-exist:nondeg_ip1}
		\liminf_{v\to\infty}\ip(v)>A
	\end{equation}
	and
	\begin{equation}\label{eq-exist:nondeg_ip2}
		\int_0^{v_0}\frac{dv}{\ip(v)}<\infty\,\ \text{for some}\ v_0>0.
	\end{equation}
	Then for any smooth precompact domain $E_0$ with $|\p E_0|<A$, there exists a weak solution $u$ of IMCF on $M$ with initial condition $E_0$, such that $E_t(u)$ is precompact for all $0\leq t\leq\log\big(A/|\p E_0|\big)$. If $E_0\subset B(x_0,r_0)$ for a radius $r_0$, then we have the quantitative bound $E_t(u)\subset B(x_0,R)$ where
	\begin{equation}\label{eq-exist:bound_main_thm}
		R=r_0+2+(2+e^t)\int_0^{\isp^{-1}(e^t|\p E_0|)}\frac{dv}{\ip(v)}<\infty.
	\end{equation}
\end{theorem}

\subsection{Proof of Theorem \ref{thm-exist:main_thm}}

For simplicity, we denote $B(r)=B(x_0,r)$ for $r\in\RR$. For each integer $k\in\NN$, $k>r_0$, let $W_k$ be a connected smooth domain with $B(k+\frac14)\subset W_k\subset B(k+\frac12)$. For each $k$, we construct a new manifold $M_k$ by smoothly attaching a metric cone $\p W_k\times[0,\infty)$ to $W_k$. The detailed construction is as follows. Choose $0<\delta\leq\frac18$ such that $2\delta$ is smaller than the normal injectivity radius of $\p W_k$. The choice of $\delta$ depends on $k$, which we make implicit for brevity. Let $\eta:[0,\infty)\to[0,\infty)$ be a smooth cutoff function such that $\eta|_{[0,1/2]}\equiv1$ and $\eta|_{[3/4,\infty)}\equiv0$. Let $g=dr^2+h(r,x)$ ($0\leq r\leq\delta$) be the metric expression in the $\delta$-collar neighborhood of $W_k$ (the positive $r$-direction pointing outward), induced by the normal exponential map. Thus $h(0,x)=g|_{\p W_k}$. Consider the new tensor
\[h'(r,x) =
(1+e^{-1/r}) \, \eta(\frac r\delta) \, h(r,x)
+ C\big(1-\eta(\frac r\delta)\big) \, \big(\frac r\delta\big)^2 \, h(\delta,x).\]
By choosing the constant $C$ sufficiently large, $h'$ satisfies the properties
\begin{equation}\label{eq-exist:properties_Mk}
	\left\{\begin{aligned}
		& \text{(1)}\ h'(r,x)>h(r,x)\text{ for all }0<r\leq\delta. \\
		& \text{(2)}\ h'(r,x)=\big(\frac r\delta\big)^2h'(\delta,x)\text{ for all }r\geq\delta.
	\end{aligned}\right.
\end{equation}
where $h'(r,x)>h(r,x)$ means that $\big(h'_{ij}(r,x)-h_{ij}(r,x)\big)dx^idx^j$ is positive-definite (where $1\leq i,j\leq n-1$). Let $M_k=W_k\cup\big(\p W_k\times[0,\infty)\big)$, endowed with the metric $g_k$ that coincides with $g$ on $W_k$ and equals to $dr^2+h'(r,x)$ on $\p W_k\times[0,\infty)$. Thus $g_k$ is a smooth metric on $M_k$. There is a smooth map
\begin{equation}\label{eq-exist:def_of_Phi}
	\Phi: W_k\cup\big(\p W_k\times[0,\delta)\big)\to M,
\end{equation}
composed of the identity map on $W_k$ and the normal exponential map on $\p W_k\times[0,\delta)$. By item (1) in \eqref{eq-exist:properties_Mk}, $\Phi$ is 1-Lipschitz (with the source metric $g_k$ and target metric $g$), and $\Phi|_{\p W_k\times(0,\delta)}$ strictly decreases the area of any hypersurface in $\p W_k\times(0,\delta)$. Finally, $(n-1)\log(r/\delta)$ is a smooth solution to the IMCF on $M_k$ with initial condition $W_k\cup\big(\p W_k\times[0,\delta)\big)$. By Theorem \ref{thm-imcf:subsol_implies_exist}, there is a proper weak solution $u_k$ to the IMCF on $M_k$, with initial condition $E_0$. Through a chain of lemmas below, we will show that $\min\big(u_k,\log(A/|\p E_0|)\big)$ is the desired weak solution for sufficiently large $k$. \\

For each $k\in\NN$, $k>r_0$, define the first escaping time of $u_k$ as follows:
\begin{equation}\label{eq-exist:def_Tk}
	T_k=\sup\big\{t\geq0: E_t(u_k)\subset B(k)\big\}.
\end{equation}
Clearly the supremum in \eqref{eq-exist:def_Tk} is achieved, so $E_{T_k}(u_k)\subset B(k)$.

Note the notational subtlety here. The precise statement for ``$E_t(u_k)\subset B(k)$'' is ``$E_t(u_k)\subset W_k$ and its identical image in $M$ is contained in $B(k)$''. For brevity of statements, we keep the simplified notation for the rest of the proof.
If a set $E\subset M_k$ is actually contained in $W_k$, then clearly $|\p^*E|_g=|\p^*E|_{g_k}$.

The following lemma controls the jumping behavior of $u_k$, and is where the properties \eqref{eq-exist:properties_Mk} of $M_k$ are used. The nice behavior of $u_k$ depends on the particular construction of $M_k$. For instance, if one creates a thin neck when attaching the exterior part to $W_k$, then $u_k$ will quickly jump to the thin neck and ignore the geometry inside $W_k$.

\begin{lemma}\label{lemma-exist:jumping_control}
	Given $r>0$ and $0\leq t\leq\log\big(A/|\p E_0|\big)$, let
	\begin{equation}\label{eq-exist:jumping_control}
		R=r+1+2\int_0^{\isp^{-1}(e^t|\p E_0|)}\frac{dv}{\ip(v)}.
	\end{equation}
	Then $R$ is finite under the isoperimetric conditions \eqref{eq-exist:nondeg_ip1} \eqref{eq-exist:nondeg_ip2}. For all $k>R$ the following statement holds: if $E_t(u_k)\subset B(r)$, then $E_t^+(u_k)\subset B(R)$.
\end{lemma}
\begin{proof}
	The finiteness of $R$ follows from Lemma \ref{lemma-sip:isp_nonzero} and \eqref{eq-exist:nondeg_ip1} \eqref{eq-exist:nondeg_ip2}. Note that $\p E_t^+(u_k)\setminus\p E_t(u_k)$ is a $g_k$-minimal surface (rigorously speaking, the support of a $g_k$-stationary integral varifold) in $M_k$, by Lemma \ref{lemma-imcf:Et_properties} (1), Definition \ref{def-hull:max_vol_sol} and Theorem \ref{def-hull:minimizing_hull}. Observe on the other hand that $\p W_k\times[\delta,\infty)$ is foliated by strictly convex hypersurfaces, hence $E_t^+(u_k)\subset W_k\cup\big(\p W_k\times(0,\delta)\big)$ by a strong maximum principle of Solomon-White \cite[Theorem 4]{White_2010}. Suppose that $E_t^+(u_k)$ has nonempty intersection with $\p W_k\times(0,\delta)$. Then the map $\Phi$ defined in \eqref{eq-exist:def_of_Phi} strictly decreases the area of $E_t^+(u_k)$. Denote $F_1=\Phi(E_t^+(u_k))$, we have by Lemma \ref{lemma-imcf:Et_properties} (4)
	\[|\p E_t(u_k)|_g=|\p E_t(u_k)|_{g_k}\geq|\p E_t^+(u_k)|_{g_k}>|\p F_1|_g.\]
	By the construction of $W_k$, we have $F_1\subset B(k+1)$. Now let $F_2$ be any minimizer of the following double obstacle problem:
	\[|\p^*F_2|_g=\inf\big\{|\p^*F|_g: E_t(u_k)\subset F\subset B(k+1)\big\}.\]
	Therefore $|\p^*F_2|_g\leq|\p F_1|_g<|\p E_t^+(u_k)|_{g_k}$. In particular, $|\p^*F_2|_g\leq e^t|\p E_0|\leq A$. Applying Lemma \ref{lemma-sip:apriori_barrier} we obtain $F_2\subset B(R)$, thus $F_2$ is a valid perimeter competitor for $E_t^+(u_k)$. This implies that $E_t^+(u_k)$ does not solve the least area problem outside $E_t(u_k)$ in $M_k$, which contradicts Lemma \ref{lemma-imcf:Et_properties} (2) and Theorem \ref{def-hull:minimizing_hull}. Hence $E_t^+(u_k)\subset W_k$, and Lemma \ref{lemma-sip:apriori_barrier} again implies $E_t^+(u_k)\subset B(R)$. This proves the lemma.
\end{proof}

\begin{cor}\label{cor-exist:Tk>0}
	There is $k_0\in\NN$ such that $T_k>0$ for all $k\geq k_0$.
\end{cor}
\begin{proof}
	Apply Lemma \ref{lemma-exist:jumping_control} with $t=0$, $r=r_0$. For each $k$ we have $E_0^+(u_k)\subset B(R)$ for all $k>R$, where $R$ is given by \eqref{eq-exist:jumping_control}. By Lemma \ref{lemma-imcf:Et_properties} (3) we have $E_\epsilon(u_k)\subset B(R+1)$ for some small $\epsilon>0$ (which may depend on $k$). This proves the lemma with $k_0=\lfloor R\rfloor+2$.
\end{proof}

For $i=1,2$, define
\begin{equation}\label{eq-exist:def_Rk}
	R_i(t)=r_0+2i+(2i+e^t)\int_0^{\isp^{-1}(e^t|\p E_0|)}\frac{dv}{\ip(v)}.
\end{equation}

Thus $R_i$ is finite when $t\leq\log\big(A/|\p E_0|\big)$, under the conditions \eqref{eq-exist:nondeg_ip1} \eqref{eq-exist:nondeg_ip2}. The following a priori diameter bound for $u_k$ lies at the heart of the proof.

\begin{lemma}\label{lemma-exist:apriori_bound_IMCF}
	We have $E_t(u_k)\subset B(R_1(t))$ for all $k\geq k_0$ and all $t$ that satisfies
	\begin{equation}\label{eq-exist:aux3}
		0<t\leq\min\big\{T_k,\log\big(A/|\p E_0|\big)\big\}.
	\end{equation}
\end{lemma}
\begin{proof}
	Fix a time $t$ satisfying \eqref{eq-exist:aux3}. We may assume $R_1(t)<k$, otherwise the statement is trivial since $E_t(u_k)\subset B(k)$. The proof here is analogous to Lemma \ref{lemma-sip:apriori_barrier}, where the inequality \eqref{eq-sip:stable2} is replaced by \eqref{eq-exist:excess_ineq} below.
	
	Suppose by contradiction that $|E_t(u_k)\setminus B(R_1(t))|>0$. Let us recall by Remark \ref{rmk-sip:level_set} that $E_{t'}=E_{t'}^{(1)}$ for all $t'>0$. For almost every $\rho\geq r_0$, we have $\H^{n-1}(\p E_t(u_k)\cap\p B(\rho))=0$ and $\p B(\rho)$ is smooth up to zero $(n-1)$-Hausdorff measure, see Remark \ref{rmk-sip:geo_sph}. For all such $\rho$ and all $t'\in[0,t]$ define
	\[V(t',\rho)=\big|E_{t'}(u_k)\setminus\bar{B(\rho)}\big|_g,\ \,
	S(t',\rho)=\big|E_{t'}(u_k)\cap\p B(\rho)\big|_g,\ \,
	A(t',\rho)=\big|\p E_{t'}(u_k)\setminus\bar{B(\rho)}\big|_g.\]
	Comparing $J_{u_k}(E_t(u_k))\leq J_{u_k}(E_t(u_k)\cap B(\rho))$, by \eqref{eq-sip:intersection_peri} and coarea formula we have
	\begin{equation}\label{eq-exist:iterated}
		A(t,\rho) \leq S(t,\rho)+\int_{E_t(u_k)\setminus\bar{B(\rho)}}|\D u_k|_g
		= S(t,\rho)+\int_0^t A(\tau,\rho)\,d\tau.
	\end{equation}
	For the same reason, we have
	\begin{equation}\label{eq-exist:aux1}
		A(t',\rho) \leq S(t',\rho)+\int_0^{t'} A(\tau,\rho)\,d\tau
	\end{equation}
	whenever $\H^{n-1}(\p E_{t'}(u_k)\cap\p B(\rho))=0$. Since \ref{eq-exist:aux1} holds for all but countably many $t'$ for each $\rho$, we can apply it in the right hand side of \eqref{eq-exist:iterated} and obtain
	\begin{equation}\label{eq-exist:iterated2}
		\begin{aligned}
			A(t,\rho) &\leq S(t,\rho)+\int_0^t\Big[S(\tau,\rho)+\int_0^\tau A(\tau',\rho)\,d\tau'\Big]\,d\tau \\
			&\leq S(t,\rho)
			+ \int_0^t S(\tau,\rho)\,d\tau
			+ \int_0^t(t-\tau)A(\tau,\rho)\,d\tau.
		\end{aligned}
	\end{equation}
	Applying \eqref{eq-exist:aux1} again to the right hand side of \eqref{eq-exist:iterated2}, we obtain
	\[A(t,\rho) \leq S(t,\rho)
	+ \int_0^t S(\tau,\rho)\,d\tau
	+ \int_0^t (t-\tau)S(\tau,\rho)\,d\tau
	+ \frac12 \int_0^t (t-\tau)^2A(\tau,\rho)\,d\tau.\]
	Iterating this process, we obtain
	\[A(t,\rho) \leq S(t,\rho)
	+ \sum_{j=0}^{m-1}\frac1{j!}\int_0^t(t-\tau)^j S(\tau,\rho)\,d\tau
	+ \frac1{m!}\int_0^t(t-\tau)^m A(\tau,\rho)\,d\tau.\]
	Since $A(\tau,\rho)\leq e^\tau|\p E_0|_g\leq e^t|\p E_0|_g$, the remainder term vanishes when $m\to\infty$. Hence
	\begin{equation}\label{eq-exist:excess_ineq}
		A(t,\rho) \leq S(t,\rho)+\int_0^t e^{t-\tau}S(\tau,\rho)\,d\tau
		\leq e^tS(t,\rho),
	\end{equation}
	where we note that $S(\tau,\rho)$ is non-decreasing in $\tau$. Next, for almost every $\rho$ we have
	\begin{equation}\label{eq-exist:ip}
		A(t,\rho)+S(t,\rho)=|\p^*(E_t(u_k)\setminus\bar{B(\rho)})|_g\geq\ip(V(t,\rho))
	\end{equation}
	by Lemma \ref{lemma-sip:intersect_decomp} and the remarks that follows. Combining \eqref{eq-exist:excess_ineq} \eqref{eq-exist:ip} and the coarea formula, we have
	\begin{equation}\label{eq-exist:ode_ineq}
		\frac d{d\rho}V(t,\rho)=-S(t,\rho)\leq-\frac1{1+e^t}\ip(V(t,\rho))
		\quad \text{for a.e.}\,\rho.
	\end{equation}
	Since $V(t,R_1(t))>0$ by our assumption, we can integrate \eqref{eq-exist:ode_ineq} and use the definition of $R_1(t)$ to obtain
	\begin{equation}\label{eq-exist:contradiction}
		\int_{V(t,R_1(t))}^{V(t,r_0)}\frac{dv}{\ip(v)} \geq
		\frac{R_1(t)-r_0}{1+e^t}>\int_0^{\isp^{-1}(e^t|\p E_0|)}\frac{dv}{\ip(v)}.
	\end{equation}
	However, we have $V(t,r_0)\leq|E_t(u_k)|\leq\isp^{-1}(|\p E_t(u_k)|)\leq\isp^{-1}(e^t|\p E_0|)$, which leads to a contradiction with \eqref{eq-exist:contradiction}.
\end{proof}

For convenience, we denote $\tilde T=\log(A/|\p E_0|)$.

\begin{cor}\label{cor-exist:lower_bound_Tk}
	There exists $k_1\in\NN$ such that $T_k\geq\tilde T$ for all $k\geq k_1$. Furthermore, $E_{\tilde T}(u_k)$ is outward minimizing in $M$ for all $k\geq k_1$.
\end{cor}
\begin{proof}
	Choose $k_1>R_2(\tilde T)$ (in particular, $k_1>k_0$). Suppose that $k\in\NN$ satisfies $T_k<\tilde T$. We will show that $k<k_1$, which proves the first statement. Applying Lemma \ref{lemma-exist:apriori_bound_IMCF} we obtain $E_{T_k}(u_k)\subset B(R_1(T_k))$. If further $R_2(T_k)<k$, then Lemma \ref{lemma-exist:jumping_control} applies to yield $E_{T_k}^+(u_k)\subset B(R_2(T_k)-1)$. Therefore $E_{T_k+\epsilon}(u_k)\subset B(R_2(T_k))\subset B(k)$ for some small $\epsilon>0$, by Lemma \ref{lemma-imcf:Et_properties} (3). This contradicts the maximality of $T_k$, hence $R_2(T_k)\geq k$. We have
	\[k\leq R_2(T_k)<R_2(\tilde T)<k_1,\]
	and the first statement follows. For each $k\geq k_1$, Lemma \ref{lemma-exist:apriori_bound_IMCF} implies $E_{\tilde T}(u_k)\subset B(R_1(\tilde T))$. By Corollary \ref{cor-sip:minimizing_hull}, $E_{\tilde T}(u_k)$ admits a strictly outward minimizing hull that is contained in $B(R_2(\tilde T)-1)$. By Lemma \ref{lemma-hull:max_vol_sol}, if $E_{\tilde T}(u_k)$ is not outward minimizing in $M$ then it is not outward minimizing in $B(R_2(\tilde T)-1)$. The latter does not hold since $u$ is a weak solution in $B(R_2(\tilde T))\subset B(k)$. This proves the second statement.
\end{proof}

\begin{proof}[Proof of Theorem \ref{thm-exist:main_thm}] {\ }
	
	Given all the settings described above, we choose $k_1$ as in Corollary \ref{cor-exist:lower_bound_Tk} and consider the function
	\[u(x)=\left\{\begin{aligned}
		& \min(u_{k_1}(x),\tilde T)\qquad(x\in B(k_1)), \\
		& \tilde T\qquad(x\notin B(k_1)).
	\end{aligned}\right.\]
	Since $E_{\tilde T}(u_{k_1})\subset B(R_1(\tilde T))\subset\subset B(k_1)$ by Lemma \ref{lemma-exist:apriori_bound_IMCF}, $u$ is a continuous function. The quantitative bound \eqref{eq-exist:bound_main_thm} for $0\leq t\leq\tilde T$ is inherited from Lemma \ref{lemma-exist:apriori_bound_IMCF}. It remains to show that $u$ is a weak solution with initial condition $E_0$. Given $0<t\leq\tilde T$ and a competitor set $E$ such that $E\Delta E_t(u)\subset K\subset\subset M\setminus E_0$. By the outward minimizing property in Corollary \ref{cor-exist:lower_bound_Tk}, we have
	\[J_u^K(E\cup E_{\tilde T}(u))\geq J_u^K(E_{\tilde T}(u)).\]
	Then by \eqref{eq-sip:cap_cup} we have
	\[J_u^K(E)\geq J_u^K(E\cap E_{\tilde T}(u))+J_u^K(E\cup E_{\tilde T}(u))-J_u^K(E_{\tilde T}(u))\geq J_u^K(E\cap E_{\tilde T}(u)).\]
	Since $u$ is a weak solution in $B(k_1)$, we have
	\[J_u^K(E\cap E_{\tilde T}(u))\geq J_u^K(E_t(u)).\]
	It follows that $J_u^K(E)\geq J_u^K(E_t(u))$, hence $u$ is a weak solution.
\end{proof}

\subsection{Proof of Theorem \ref{thm-intro:existence}}

Assume the conditions of Theorem \ref{thm-intro:existence}. For each $l\in\NN_+$, we apply Theorem \ref{thm-exist:main_thm} with the choice $A=e^l|\p E_0|$. We obtain a sequence of weak solutions $u^l$ with initial condition $E_0$, such that $E_l(u^l)$ is precompact in $M$, and the quantitative diameter bound \eqref{eq-exist:bound_main_thm} holds for $u^l$ whenever $t\leq l$. By Theorem \ref{thm-imcf:max_principle} (2), for two integers $l<l'$ we have $u^l=u^{l'}$ on $E_l(u^l)$. Therefore, the function
\[u(x)=\lim_{l\to\infty}u^l(x)\]
is defined on $\bigcup_{l\in\NN}E_l(u^l)$ (which is $M$, by the local gradient estimate \eqref{eq-imcf:gradient_estimate}). For each $t>0$ we have $E_t(u)=E_t(u^{\lceil t\rceil})$, hence \eqref{eq-exist:bound_main_thm} holds for $u$ as well. In particular, $u$ is proper. It remains to show that $u$ is a weak solution. Consider any energy competitor $E$ with $E\Delta E_t(u)\subset K\subset\subset M\setminus\bar{E_0}$. Since $K$ is compact, it is covered by finitely many domains $E_l(u^l)$, hence is included in a single domain $E_l(u^l)$ for some $l$ (since the sets $E_l(u^l)$ is increasing). By construction we have $u=u^l$ in $E_l(u^l)$. Using the fact that $u^l$ is a weak solution, we have
\[J_u^K(E)=J_{u^l}^K(E)\geq J_{u^l}^K(E_t(u^l))=J_u^K(E_t(u)).\]
Therefore $u$ is a weak solution. This proves existence, and uniqueness follows from the maximum principle. \hfill$\Box$


\vspace{12pt}

\noindent\textit{{Kai Xu,}}

\vspace{2pt}

\noindent\textit{{Department of Mathematics, Duke University, Durham, NC 27708, USA,}}

\vspace{2pt}

\noindent\textit{Email address: }\href{mailto:kx35@math.duke.edu}{kx35@math.duke.edu}.

\end{document}